\documentclass{article}
\usepackage{mathrsfs}
\usepackage{amsmath}
\usepackage{amssymb, amsmath}
\usepackage{graphicx}
\catcode`\@=11 \@addtoreset{equation}{section}

\catcode`\@=12
\usepackage{color}

\newtheorem{Theorem}{Theorem}[section]
\newtheorem{Proposition}{Proposition}[section]
\newtheorem{Lemma}{Lemma}[section]
\newtheorem{Corollary}{Corollary}[section]
\newtheorem{Definition}{Definition}[section]

\newtheorem{Remark}{Remark}[section]

\newcommand{\newcom}{\newcommand}
\newcommand{\bTheorem}[1]{
\begin{Theorem} \label{T#1} }
\newcommand{\eT}{\end{Theorem}}

\newcommand{\bProposition}[1]{
\begin{Proposition} \label{P#1}}
\newcommand{\eP}{\end{Proposition}}

\newcommand{\bLemma}[1]{
\begin{Lemma} \label{L#1} }
\newcommand{\eL}{\end{Lemma}}

\newcommand{\bCorollary}[1]{
\begin{Corollary} \label{C#1} }
\newcommand{\eC}{\end{Corollary}}

\newcommand{\beq}{\begin{equation}}
\newcommand{\eeq}{\end{equation}}
\newcom{\ben}{\begin{eqnarray}}
\newcom{\een}{\end{eqnarray}}
\newcom{\beno}{\begin{eqnarray*}}
\newcom{\eeno}{\end{eqnarray*}}
\newcom{\bali}{\begin{aligned}}
\newcom{\eali}{\end{aligned}}

\newcommand{\bFormula}[1]{
\begin{equation} \label{#1}}
\newcommand{\eF}{\end{equation}}

\newcommand{\f}{\frac}
\newcommand{\Om}{\Omega}

\newcommand{\p}{\partial}

\newcommand{\vr}{\varrho}

\newcommand{\vt}{\vartheta}
\newcommand{\vu}{\vc{u}}

\newcommand{\vv}{\vc{v}}

\newcommand{\vc}[1]{{\boldsymbol #1}}
\newcommand{\Div}{{\rm div}}
\newcommand{\Grad}{\nabla}

\newcommand{\dx}{{\rm d} x}

\newcommand{\dt}{{\rm d} t }
\newcommand{\ds}{{\rm d} s}

\newcommand{\dxdt}{\dx\dt}

\newcommand{\ep}{\varepsilon}

\font\F=msbm10 scaled 1000
\newcommand{\R}{\mbox{\F R}}

\newcommand\Cbox[2]{%
    \newbox\contentbox%
    \newbox\bkgdbox%
    \setbox\contentbox\hbox to \hsize{%
        \vtop{
            \kern\columnsep
            \hbox to \hsize{%
                \kern\columnsep%
                \advance\hsize by -2\columnsep%
                \setlength{\textwidth}{\hsize}%
                \vbox{
                    \parskip=\baselineskip
                    \parindent=0bp
                    #2
                }%
                \kern\columnsep%
            }%
            \kern\columnsep%
        }%
    }%
    \setbox\bkgdbox\vbox{
        \color{#1}
        \hrule width  \wd\contentbox %
               height \ht\contentbox %
               depth  \dp\contentbox
        \color{black}
    }%
    \wd\bkgdbox=0bp%
    \vbox{\hbox to \hsize{\box\bkgdbox\box\contentbox}}%
    \vskip\baselineskip%
}

\begin{document}


\title{\bf Global weak solutions to a two-dimensional compressible MHD equations of viscous non-resistive fluids}

\author{
Yang Li \\ Department of Mathematics, \\ Nanjing University, Nanjing 210093, China \\ lymath@smail.nju.edu.cn\\
Yongzhong Sun \\ Department of Mathematics, \\ Nanjing University, Nanjing 210093, China \\ sunyz@nju.edu.cn
}

\maketitle
{ \centerline {\bf Abstract}}

{We consider a two-dimensional MHD model describing the evolution of viscous, compressible and electrically conducting fluids under the action of vertical magnetic field without resistivity. Existence of global weak solutions is established for any adiabatic exponent $\gamma \geq 1$. Inspired by the approximate scheme proposed in \cite{FNP}, we consider a two-level approximate system with artificial diffusion and pressure term. At the first level, we prove global well-posedness of the regularized system and establish uniform-in-$\ep$ estimates to the regular solutions. At the second level, we show global existence of weak solutions to the system with artificial pressure by sending $\ep $ to $0$ and deriving uniform-in-$\delta$ estimates. Then global weak solution to the original system is constructed by vanishing $\delta$. The key issue in the limit passage is the strong convergence of approximate sequence of the density and magnetic field. This is accomplished by following the technique developed in \cite{FNP,LS} and using the new technique of variable reduction developed by Vasseur et al. \cite{VWY} in the context of compressible two-fluid model so as to handle the cross terms.}

{\bf Keywords: }{Non-resistive MHD equations; Weak solutions; Global existence}

{\bf Mathematics Subject Classification: }{35M10, 35G31, 35D30 }

\section{Introduction}
The motion of electrically conducting fluids under the interactions of the magnetic field is described by the equations of  magnetohydrodynamics (MHD).  In Eulerian coordinates, a simplified and well-accepted model for the three-dimensional compressible MHD equations reads as follows (see \cite{HC}):
\beq\label{p1}
\p_t\vr + \Div(\vr\vv)=0,
\eeq
\beq\label{p2}
\p_t(\vr\vv)+\Div(\vr\vv \otimes \vv)+\Grad{p}=\mu\Delta\vv + (\mu +\lambda)\Grad{\Div \vv}
+(\Grad \times \vc{B})\times \vc{B},
\eeq
\beq\label{p3}
\p_t\vc{B}=\Grad \times (\vv \times \vc{B})-\nu \Grad \times (\Grad \times \vc{B}),
\eeq
\beq\label{p4}
\Div \vc{B}=0.
\eeq
Here the unknowns $\vr,\,\vv \in \R^3,\,p$ and $\vc{B}\in\R^3$ denote the density of fluids, the velocity field, the pressure and the magnetic field, respectively. The viscosity coefficients $\mu$ and $\lambda$ satisfy
\[
\mu>0 ,\,\,\lambda+2\mu > 0.
\]
Moreover, $\nu\ge 0$ is the resistivity coefficient which represents magnetic diffusion of the magnetic field. The compressible flow is assumed to be isentropic (isothermal). This means the pressure $p$ is prescribed through the following constitutive relation:
\beq\label{p5}
p(\vr)=a\vr^{\gamma},
\eeq
where $a$ is a positive constant and the adiabatic exponent $\gamma \geq1$.

By setting $\vc{B}=\vc{0}$, the equations (\ref{p1})-(\ref{p5}) reduce to the isentropic (isothermal) compressible Navier-Stokes system. Global existence of weak solutions with smallness restrictions imposed on the initial data was proved by Hoff \cite{HD2,HD1}; while global existence for spherically symmetric solutions was obtained by Jiang and Zhang \cite{JZP}. Global existence of weak solutions with large data has been verified by Lions \cite{LS} and refined by Feireisl et al. \cite{FE1,FNP}. We refer to \cite{FE2,FN1} for the global existence theory of weak solutions to the Navier-Stokes-Fourier system. The construction of global weak solutions becomes much more complicated when the mutual interactions between the electrically conducting fluids and the magnetic field are taken into account. Global existence of weak solutions to the three-dimensional isentropic compressible MHD system was proved by Hu and Wang \cite{HW}. Ducomet and Feireisl \cite{DF} considered the heat-conducting fluids together with the influence of radiation, and obtained the global existence of weak solutions with finite energy initial data.

Remarkably, the results in \cite{DF,HW} are based on the assumption that the resistivity coefficient $\nu$ is positive. In practical models such as plasmas in general $\nu$ is very small. Therefore, it is natural to consider the compressible non-resistive MHD equations. In fact, the lack of diffusion mechanism for the magnetic field brings extra obstacle in constructing global solutions. Even in the one-dimensional case, global existence and uniqueness of smooth solutions with large data has recently been obtained by Jiang and Zhang \cite{JZW}; see also \cite{LY1} for the treatment of planar non-resistive MHD equations. Wu and Wu \cite{WW} showed the global solvability of small smooth solutions to the two-dimensional compressible non-resistive MHD equations; the verification in three space dimensions can be found in Jiang and Jiang \cite{JJ1,JJ2}. For investigations on the incompressible non-resistive MHD system, we refer to \cite{LXZ,RWZ,TW,XZ,ZT} and the references therein. As far as the weak solution is concerned, the only global existence result is restricted to the one-dimensional regime \cite{LYS}. The existence of global weak solutions to the multi-dimensional compressible non-resistive MHD equations with large data is largely open. The purpose of this paper is to move a step in this direction. Due to mathematical challenges, we consider the case that the motion of fluids takes place in the plane $\R^2$, and the magnetic field acts on fluids only in the vertical direction. Precisely, by choosing
\[
\vv=(u^1(x_1,x_2,t),u^2(x_1,x_2,t),0)=:(\vu,0)
\]
and
\[
\vc{B}=(0,0,b(x_1,x_2,t)),\,\,\,\vr=\vr(x_1,x_2,t),
\]
the equations (\ref{p1})-(\ref{p5}) (with $\nu =0$) read as follows:
\beq\label{p6}
\p_t\vr + \Div(\vr\vu)=0,
\eeq
\beq\label{p7}
\p_t(\vr\vu)+\Div(\vr\vu \otimes \vu)+\Grad{\left(a\vr^{\gamma}+\f{1}{2}b^2\right)}=\mu\Delta\vu + (\mu +\lambda)\Grad{\Div \vu},
\eeq
\beq\label{p8}
\p_t b+ \Div(b\vu)=0.
\eeq

The equations (\ref{p6})-(\ref{p8}) are reminiscent of the compressible two-fluid model. In three space dimensions, global existence of weak solutions was first obtained by Vasseur et al. \cite{VWY} when the pressure satisfies $\gamma$-laws (a very recent treatment can be found in \cite{Wen}). We refer to Bresch et al. \cite{BMZ} in the context of semi-stationary Stokes system with algebraic pressure closure; see also Novotn\'{y} and Pokorn\'{y} \cite{NP} for the non-stationary case with general form of pressure. Yao et al. \cite{YZZ} obtained the existence and  asymptotic behavior of global weak solutions to the two-dimensional liquid-gas model with smallness restrictions imposed on the initial data. Also we refer to \cite{MMM} for the existence of global weak solutions to the compressible Navier-Stokes equations with entropy transport, which partly motivated the present work. Concerning the one-dimensional two-fluid model, we refer to \cite{EK1,EWZ1} and the references therein.

To fix ideas, we assume the motion of fluids takes place in a bounded regular domain $\Om \subset \R^2$ and the adiabatic exponent $\gamma>1$. The equations (\ref{p6})-(\ref{p8}) are supplemented with the initial conditions:
\beq\label{p9}
\vr(x,0)=\vr_0(x),\,\, \vu (x,0)=\vu_0 (x),\,\,b(x,0)=b_0(x), \,\,\,x\in \Om,
\eeq
together with the no-slip boundary condition:
\beq\label{p10}
\vu|_{\p\Om}=\vc{0}.
\eeq

Formally, taking inner product of (\ref{p7}) with $\vu$ and integrating by parts over $\Om$, with the help of (\ref{p6}) and (\ref{p8}), gives the energy inequality
\beq\label{p11}
\f{d}{dt} \mathcal{E}(t)+ \int_{\Om}\left( \mu |\nabla \vu|^2+(\mu+\lambda)(\Div \vu)^2 \right)\dx \leq 0,
\eeq
where the total energy $\mathcal{E}(t)$ is given by
\[
\mathcal{E}(t):= \int_{\Om}\left( \f{1}{2}\vr |\vu|^2 +\f{a}{\gamma-1}\vr ^{\gamma}+\f{1}{2}b^2
\right)(x,t) \dx.
\]

Now we introduce the definition of weak solutions.
\begin{Definition}\label{din}
A triple $(\vr,\vu,b)$ is said to be a weak solution to (\ref{p6})-(\ref{p10}) on $\Om\times (0,T)$ provided that
\begin{itemize}
\item {  $\vr \geq 0, \,\,\vr \in L^{\infty}(0,T;L^{\gamma}(\Om)) ,\,\,b\in L^{\infty}(0,T;L^2(\Om)),\,\,\vu \in L^2(0,T;H^1_0(\Om))  $,  }

\item {$(\vr,\vu,b)$ solves (\ref{p6})-(\ref{p8}) in $\mathcal{D}'(\Om\times (0,T))$. In addition, (\ref{p6}) and (\ref{p8}) are satisfied in $\mathcal{D}'(\R^2 \times (0,T))$ if $\vr,\vu$ and $b$ are set to be zero outside $\Om$,    }

\item {the energy inequality (\ref{p11}) holds in $\mathcal{D}'((0,T))$,  }

\item {(\ref{p6}) and (\ref{p8}) are satisfied in the sense of renormalized solutions, i.e.,
\[
\p_t h(\vr)+\Div(h(\vr) \vu)+(h'(\vr)\vr-h(\vr))\Div\vu=0  \,\, \text{ in } \mathcal{D}'(\Om\times(0,T)),
\]
\[
\p_t h(b)+\Div(h(b) \vu)+(h'(b)b-h(b))\Div\vu=0  \,\, \text{ in } \mathcal{D}'(\Om\times(0,T)),
\]
for any $h\in C^1(\R)$ with $h'(z)=0$ for sufficiently large value of $z$.
  }

\end{itemize}
\end{Definition}

It is a routine matter to check, by (\ref{p6})-(\ref{p8}), that any weak solution $(\vr,\vu,b)$ belongs to the regularity class
\[
\vr\in C([0,T];L^{\gamma}_{w}(\Om)),\,\,\vr \vu \in C([0,T];L^{\f{2\gamma}{\gamma+1}}_{w}(\Om)),\,\, b\in C([0,T];L^{2}_{w}(\Om)).
\]
Here, we denote by $C([0,T];L^p_{w}(\Om))$ the space of continuous functions from $[0,T]$ to $L^p(\Om)$ equipped with the weak topology.

Our main result of this paper reads as follows.
\begin{Theorem}\label{ls}
Assume that $\gamma > 1$ and $\Om\subset\R^2$ be a bounded domain of class $C^{2,\beta}$, $\beta>0$. Let  the initial data $(\vr_0,\vu_0,b_0)$ be subject to
\beq\label{p12}
0<m\leq \vr_0,b_0\leq M<\infty,
\eeq
\beq\label{p13}
\vu_0 \in L^2(\Om).
\eeq

Then there exists a global weak solution to (\ref{p6})-(\ref{p10}) emanating from $(\vr_0,\vu_0,b_0)$.
\end{Theorem}

\begin{Remark}
As a matter of fact, the same conclusion of Theorem \ref{ls} holds provided that the initial data are prescribed through
\beq\label{p14}
\vr_0 \geq 0,\,\,\,\underline{C} \,\vr_0 \leq b_0 \leq \overline{C}\,\vr_0\,\, \text{  a.e. in }\Om,
\eeq
\beq\label{p15}
\vr_0 \in L^{\gamma}(\Om),\,\, b_0 \in  L^2(\Om),\, \, \f{ \vc{m}_0}{\sqrt{\vr_0}} \in L^2(\Om),
\eeq
\beq\label{p16}
 \vc{m}_0 =\vc{0}\,\text{ a.e. in }    \{x\in \Om: \vr_0(x)=0\},
\eeq
for some constants $0< \underline{C} <\overline{C}<\infty$. Notice that in fact $(\vr_0, b_0) \in L^{\max\{2,\gamma\}}(\Om)$ according to (\ref{p14}). It turns out that (\ref{p14}) appears naturally in verifying global well-posedness of the first level approximate system due to the specific mathematical structure of our MHD model; see Section \ref{sen2}.
\end{Remark}

\begin{Remark}
With slight modifications in the proof, the conclusion of Theorem \ref{ls} remains valid if the adiabatic exponent $\gamma=1$. To see this, it suffices to observe that the square integrability of the magnetic field will guarantee that of the density, based on the assumption (\ref{p14}); see (\ref{pq9}), (\ref{ps13}), (\ref{ps24}), (\ref{pt11}).
\end{Remark}

\begin{Remark}
Basically, global weak solutions to (\ref{p6})-(\ref{p10}) are constructed by considering an approximate system with artificial viscosity terms in the continuity equations and pressure term in spirit of \cite{FNP}. The essential part is to verify the pointwise convergence of approximate densities and magnetic fields. This is a delicate issue since now the pressure term depends on two variables each satisfying the continuity equation. We appeal to the new method of variable reduction proposed by Vasseur et al. \cite{VWY} in the context of compressible two-fluid model to handle the cross terms.
\end{Remark}

The rest of this paper is organized as follows. In Section \ref{sen2}, we establish global well-posedness to the first level approximate system for (\ref{p6})-(\ref{p10}) in light of Brenner's model, which may be of independent interest. In Section \ref{sen3}, we give uniform-in-$\ep$ estimates to the strong solutions $\{(\vr_{\ep},\vu_{\ep},b_{\ep})\}_{\ep>0}$ obtained in Section \ref{sen2} and pass to the limit as $\ep \rightarrow 0$. Finally, by establishing uniform-in-$\delta$ estimates, we pass to the limit for approximate solutions $\{(\vr_{\delta},\vu_{\delta},b_{\delta})\}_{\delta>0}$ as $\delta \rightarrow 0$ to recover global weak solution to the original problem (\ref{p6})-(\ref{p10}); see Section \ref{sen4}.

\section{Global well-posedness to the first level approximate system}\label{sen2}
Inspired by the model proposed by Brenner \cite{BR} and the approximate scheme for constructing weak solutions to the compressible Navier-Stokes equations \cite{FNP}, we consider the following approximate system associated with (\ref{p6})-(\ref{p8}):
\beq\label{pq1}
\p_t\vr + \Div(\vr\vu)=\ep \Delta \vr,
\eeq
\[
\p_t(\vr\vu)+\Div(\vr\vu \otimes \vu)+\Grad{\left(a\vr^{\gamma}+\f{1}{2}b^2\right)}
+\ep \nabla \vr \cdot \nabla \vu +\delta\nabla (\vr+b)^{\Gamma}
\]
\beq\label{pq2}
=\mu\Delta\vu + (\mu +\lambda)\Grad{\Div \vu},
\eeq
\beq\label{pq3}
\p_t b+ \Div(b\vu)=\ep \Delta b,
\eeq
where $\ep$ and $\delta$ are small positive parameters; $\Gamma>1$ is a suitably large parameter. We impose homogeneous Neumann boundary condition for the density and magnetic field, together with the no-slip boundary condition for the velocity as follows:
\beq\label{pq4}
\nabla \vr \cdot \vc{n}|_{\p\Om}=\nabla b \cdot \vc{n}|_{\p\Om}=0,\,\,\vu|_{\p\Om}=\vc{0},
\eeq
where $\vc{n}$ is the unit outward normal on $\p\Om$. Moreover, the initial conditions are prescribed through
\beq\label{pq5}
\vr(x,0)=\vr_0(x),\,\,\vu (x,0)=\vc{u}_0 (x),\,\,b(x,0)=b_0(x),\,\,\,x\in \Om.
\eeq

Our main goal in this section is to show
\begin{Proposition}\label{prop}
Let $\Om\subset \R^2$ be a bounded smooth domain and $\gamma>1$. Assume that
\beq\label{pq6}
\vr_0 >0 \,\,\text{  a.e. in } \Om,\,\, \vr_0,\vr_0^{-1}\in L^{\infty}(\Om),
\eeq
\beq\label{pq7}
b_0 >0 \,\,\text{  a.e. in } \Om,\,\, b_0,b_0^{-1}\in L^{\infty}(\Om),
\eeq
\beq\label{pq8}
(\vr_0,b_0) \in H^1(\Om),\,\, \, \vu_0\in H^1_0(\Om).
\eeq

Then the initial-boundary value problem (\ref{pq1})-(\ref{pq5}) admits a unique global strong solution $(\vr,\vu,b)$ such that for any fixed $0<T<\infty$,
\beq\label{pq9}
C_{\star}\vr\leq b\leq C^{\star}\vr\, \, \text{  a.e. in }\Om\times (0,T),
\eeq
\beq\label{pq10}
\|(\vr,\vr^{-1})\|_{L^{\infty}(\Om\times (0,T))} +\|(\vr,\vu,b)\|_{L^{\infty}(0,T;H^1(\Om))}
+\|(\vr,\vu,b)\|_{L^2(0,T;H^2(\Om))}\leq \Lambda,
\eeq
where $C_{\star}= \inf\limits_{x\in\Om}\frac{b_0(x)}{\vr_0(x)}, \,C^{\star}=\sup\limits_{x\in\Om}\frac{b_0(x)}{\vr_0(x)}$.
\end{Proposition}
In this section we use $\Lambda$ to denote a positive constant which changes from line to line and depends only on $T$, the parameters $\ep,\delta,\Gamma,\mu,\lambda,a,\gamma$ appeared in (\ref{pq1})-(\ref{pq3}) and the initial data.
\begin{Remark}
In the absence of the magnetic field, the equations (\ref{pq1})-(\ref{pq3}) reduce to
\[
\p_t\vr + \Div(\vr\vu)=\ep \Delta \vr,
\]
\[
\p_t(\vr\vu)+\Div(\vr\vu \otimes \vu)+\Grad{\left(a\vr^{\gamma} \right) }
+\ep \nabla \vr \cdot \nabla \vu+\delta\nabla \vr^{\Gamma}
=\mu\Delta\vu + (\mu +\lambda)\Grad{\Div \vu}.
\]
Global well-posedness to the above system has been obtained by Cai et al. \cite{CCS}. Our MHD system (\ref{pq1})-(\ref{pq3}) is more complicated due to the mutual interactions between electrically conducting fluids and the magnetic field.
\end{Remark}

\begin{Remark}
In fact, one could establish the existence of global weak solutions to the first level approximate system (\ref{pq1})-(\ref{pq3}) by the same approach as \cite{FNP,NP,VWY} in the context of three space dimensions.

\end{Remark}

The local-in-time existence and uniqueness of strong solution to the initial-boundary value problem (\ref{pq1})-(\ref{pq5}) can be established in a standard way based on linearization and fixed point argument; see for instance \cite{NS}. Also uniqueness of global strong solutions is justified in a routine manner based on the stability result. We refer to \cite{CCS,NS} for similar argument. Consequently, to prove Proposition \ref{prop} one only needs to derive sufficient global a priori estimates.

\subsection{Global a priori estimates}\label{sen21}
Given any $0<T<\infty$, we assume that $(\vr,\vu,b)$ is a smooth solution to the initial-boundary value problem (\ref{pq1})-(\ref{pq5}) on the space-time domain $\Om\times(0,T)$ with positive smooth initial data. We split the derivation of a priori estimates into several steps.

\noindent{ \emph{Step 1. Energy estimates.}}

By taking inner product of (\ref{pq2}) with $\vu$ and making use of the continuity equations (\ref{pq1}) and (\ref{pq3}),
\[
\f{d}{dt}\int_{\Om} \left(  \f{1}{2}\vr |\vu|^2 +\f{a}{\gamma-1} \vr^{\gamma}+\f{1}{2}b^2
+\f{\delta}{\Gamma-1}(\vr+b)^{\Gamma} \right) \dx
\]
\[
+ \int_{\Om} \left( \ep a \gamma \vr^{\gamma-2} |\nabla \vr|^2 +\ep |\nabla b|^2
+\ep \delta \Gamma (\vr+b)^{\Gamma-2} |\nabla (\vr+b)|^2  \right) \dx
\]
\beq\label{pq11}
+\int_{\Om} \left(\mu |\nabla \vu|^2+ (\mu+\lambda)(\Div \vu)^2 \right)\dx =0.
\eeq
It follows that
\[
\|\vr\|_{L^{\infty}(0,T;L^{\gamma}(\Om))}+ \|b\|_{L^{\infty}(0,T;L^2(\Om))}+
\|\vr+b\|_{L^{\infty}(0,T;L^{\Gamma}(\Om))}
\]
\[
+\|\nabla \vu\|_{L^2(0,T;L^2(\Om))}+ \|\nabla b\|_{L^2(0,T;L^2(\Om))} +\|\sqrt{\vr} \vu\|_{L^{\infty}(0,T;L^2(\Om))}
\]
\beq\label{pq12}
+\|\nabla \vr^{\f{\gamma}{2}}\|_{L^2(0,T;L^2(\Om))} +
\|\nabla (\vr+b)^{\f{\Gamma}{2}}\|_{L^2(0,T;L^2(\Om))}\leq \Lambda.
\eeq

\noindent{\emph{ Step 2. Higher integrability of the density.}}

One deduces from (\ref{pq1}) that
\[
\f{d}{dt}\int_{\Om}\vr ^{q} \dx +\ep q(q-1)\int_{\Om} \vr ^{q-2}|\nabla \vr|^2 \dx
=(1-q)\int_{\Om} \vr^{q} \Div \vu \dx,
\]
for any $q\in [2,\infty)$. Using Cauchy-Schwarz's inequality and the interpolation inequality,
\[
 \left| \int_{\Om} \vr^{q} \Div \vu \dx \right|
 \leq \|\vr^q\|_{L^2(\Om)} \|\nabla \vu\|_{L^2(\Om)}
 \leq \|\vr^{ \f{q}{2} }\|_{L^4(\Om)}^2  \|\nabla \vu\|_{L^2(\Om)}
\]
\[
\leq \Lambda(q)  \|\vr^{ \f{q}{2} }\|_{L^2(\Om)} \left( \|\vr^{ \f{q}{2} }\|_{L^2(\Om)}+\| \nabla \vr^{ \f{q}{2} }\|_{L^2(\Om)} \right )\|\nabla \vu\|_{L^2(\Om)}.
\]
Hence,
\[
\f{d}{dt}\int_{\Om}\vr ^{q} \dx +\ep \int_{\Om} \vr ^{q-2}|\nabla \vr|^2 \dx
\leq \Lambda(q) \left(   \|\nabla \vu\|_{L^2(\Om)}+ \|\nabla \vu\|_{L^2(\Om)}^2   \right)
\int_{\Om}\vr ^{q} \dx.
\]
Application of Gronwall's inequality gives
\beq\label{pq13}
\|\vr\|_{L^{\infty}(0,T;L^{q}(\Om))}+\|\nabla \vr^{\f{q}{2}}\|_{L^2(0,T;L^2(\Om))}
\leq \Lambda (q),
\eeq
where we have used
\[
\|\nabla \vu\|_{L^2(0,T;L^2(\Om))} \leq \Lambda.
\]

\noindent{\emph{ Step 3. Maximum principle.}}

From (\ref{pq6})-(\ref{pq7}),
\[
 C_{\star}\vr_0\leq b_0\leq C^{\star}\vr_0.
\]
A straightforward calculation shows that $b-C_{\star}\vr$ verifies the initial-boundary value problem
\[
\p_t\left(b-C_{\star}\vr \right) +\Div \left(  \left(b-C_{\star}\vr \right) \vu    \right)
=\ep \Delta \left(b-C_{\star}\vr \right),
\]
\[
\left(b-C_{\star}\vr \right)(x,0)=b_0(x)-C_{\star}\vr_0(x),
\]
\[
\nabla \left(b-C_{\star}\vr \right)\cdot \vc{n} |_{\p \Om}=0.
\]
Therefore, we conclude from the maximum principle for parabolic equation that
\[
b \geq C_{\star}\vr  \,\,\, \text{ in } \Om \times (0,T).
\]
Similarly, it holds that
\[
b \leq C^{\star}\vr  \,\,\, \text{ in } \Om \times (0,T).
\]
In conclusion, we have obtained
\beq\label{pq14}
C_{\star}\vr\leq b\leq C^{\star}\vr \,\, \,\text{  in }\Om\times (0,T).
\eeq

\noindent{\emph{ Step 4. Higher integrability of the velocity.}}

By taking inner product of (\ref{pq2}) with $r |\vu|^{r-2}\vu$ for some $r>2$ to be determined later and integrating the resulting relation in $\Om$, one sees after integration by parts and using (\ref{pq1}) that
\[
\f{d}{dt}\int_{\Om} \vr |\vu|^r\dx + \int_{\Om} r|\vu|^{r-2} \left(\mu |\nabla \vu|^2+(\mu+\lambda)
(\Div \vu)^2 \right)\dx
\]
\[
+r(r-2)\mu \int_{\Om}|\vu|^{r-2} |\nabla |\vu||^2\dx
+r(r-2)(\mu+\lambda)\int_{\Om} \Div \vu |\vu|^{r-3} \vu \cdot \nabla |\vu| \dx
\]
\beq\label{pq15}
=r\int_{\Om} \left( a\vr^{\gamma}+\f{1}{2}b^2 + \delta (\vr+b)^{\Gamma} \right)\Div \left(|\vu|^{r-2}\vu\right)\dx.
\eeq
Obviously, by choosing
\[
r=2+\f{\mu}{8(\mu+\lambda)},
\]
it holds that
\beq\label{pq16}
\left|r(r-2)(\mu+\lambda)\int_{\Om} \Div \vu |\vu|^{r-3} \vu \cdot \nabla |\vu| \dx\right|
\leq \f{r\mu}{8} \int_{\Om} |\vu|^{r-2} |\nabla \vu|^2 \dx.
\eeq
By virtue of H\"{o}lder's inequality, the right-hand side of (\ref{pq15}) can be estimated through
\[
\left|\int_{\Om} \left( a\vr^{\gamma}+\f{1}{2}b^2 + \delta (\vr+b)^{\Gamma} \right)\Div \left(|\vu|^{r-2}\vu\right)\dx \right|
\]
\[
\leq \Lambda  \int_{\Om} \left( \vr^{\gamma}+b^2 + (\vr+b)^{\Gamma} \right) |\vu|^{r-2}|\nabla \vu| \dx
\]
\[
\leq \nu_1 \int_{\Om} |\vu|^{r-2}|\nabla \vu|^2\dx+
\Lambda(\nu_1) \| \vr^{ \gamma-\f{r-2}{2r}  }  \|_{L^r(\Om)}^2 \|\vr |\vu|^r\|_{L^1(\Om)}^{ \f{r-2}{r}   }
\]
\[
+\nu_2 \int_{\Om} |\vu|^{r-2}|\nabla \vu|^2\dx+
\Lambda(\nu_2) \| b^{ 2-\f{r-2}{2r}  }  \|_{L^r(\Om)}^2 \|b|\vu|^r\|_{L^1(\Om)}^{ \f{r-2}{r}   }
\]
\[
+\nu_3 \int_{\Om} |\vu|^{r-2}|\nabla \vu|^2\dx+
\Lambda(\nu_3) \| (\vr+b)^{ \Gamma-\f{r-2}{2r}  }  \|_{L^r(\Om)}^2 \|(\vr+b)|\vu|^r\|_{L^1(\Om)}^{ \f{r-2}{r}   }
\]
\[
\leq (\nu_1+\nu_2+\nu_3) \int_{\Om} |\vu|^{r-2}|\nabla \vu|^2\dx+
\Lambda(\nu_1) \| \vr\|_{L^{ r(\gamma-\f{r-2}{2r}) }(\Om)}^{2(\gamma-\f{r-2}{2r}) }
\|\vr |\vu|^r\|_{L^1(\Om)}^{ \f{r-2}{r}   }
\]
\[
+\Lambda(\nu_2) \| \vr\|_{L^{ r(2-\f{r-2}{2r}) }(\Om)}^{2(2-\f{r-2}{2r}) }
\|\vr |\vu|^r\|_{L^1(\Om)}^{ \f{r-2}{r}   }
+
\Lambda(\nu_3) \| \vr\|_{L^{ r(\Gamma-\f{r-2}{2r}) }(\Om)}^{2(\Gamma-\f{r-2}{2r}) }
\|\vr |\vu|^r\|_{L^1(\Om)}^{ \f{r-2}{r}   },
\]
where we have invoked (\ref{pq14}) in the third inequality. Consequently, we choose $\nu_i=\f{\mu}{6}(i=1,2,3)$ and make use of (\ref{pq13}) to find that
\[
\left|\int_{\Om} \left( a\vr^{\gamma}+\f{1}{2}b^2 + \delta (\vr+b)^{\Gamma} \right)\Div \left(|\vu|^{r-2}\vu\right)\dx \right|
\]
\beq\label{pq17}
\leq \f{\mu}{2}\int_{\Om} |\vu|^{r-2}|\nabla \vu|^2\dx+
\Lambda \|\vr |\vu|^r\|_{L^1(\Om)}^{ \f{r-2}{r}   }.
\eeq
Combining (\ref{pq16})-(\ref{pq17}), one infers from (\ref{pq15}) that
\beq\label{pq18}
\|\vr^{\f{1}{r} }\vu \| _{L^{\infty}(0,T;L^{r}(\Om))} + \|  \nabla |\vu|^{ \f{r}{2} } \|_{L^{2}(0,T;L^{2}(\Om))} \leq \Lambda .
\eeq
As a direct consequence of Sobolev's embedding inequality, we finally arrive at
\beq\label{pq19}
\|\vu\| _{L^{r}(0,T;L^{ \f{rq}{2} }(\Om))} \leq \Lambda,
\eeq
for any $q\in [1,\infty)$.

\noindent{\emph{ Step 5. Pointwise bounds and higher order estimates for the density and magnetic field.}}

One may choose $q$ in (\ref{pq19}) suitably large such that
\[
\f{1}{q}< \f{r}{4} \left(1-\f{2}{r}\right),
\]
i.e.,
\[
\f{2}{r}+\f{4}{rq} <1.
\]
By applying Propositions 2.1-2.2 in \cite{CCS} to the parabolic equation (\ref{pq1}), it follows that
\beq\label{pq20}
\|(\vr,\vr^{-1})\|_{L^{\infty}(\Om\times (0,T))} \leq \Lambda.
\eeq
Obviously, (\ref{pq14}) and (\ref{pq20}) ensure that
\beq\label{pq21}
\|(b,b^{-1})\|_{L^{\infty}(\Om\times (0,T))} \leq \Lambda.
\eeq
Now we can obtain higher order energy estimates for $(\vr,b)$. To do this, we conclude from (\ref{pq12}) and (\ref{pq20}) that
\[
\|\vu\|_{L^{\infty}(0,T;L^{2}(\Om))} \leq \Lambda,
\]
which, together with
\[
\|\nabla \vu\|_{L^2(0,T;L^2(\Om))} \leq \Lambda,
\]
gives rise to
\beq\label{pq22}
\|\vu\|_{L^4(0,T;L^4(\Om))} \leq \Lambda.
\eeq
Multiplying (\ref{pq1}) by $\p_t\vr$ and integrating by parts over $\Om$ shows that
\[
\f{\ep}{2} \f{d}{dt}\int_{\Om}|\nabla \vr|^2 \dx +\int_{\Om}(\p_t\vr)^2\dx
=-\int_{\Om}\Div(\vr \vu)\p_t \vr \dx
\]
\[
\leq \f{1}{2}\int_{\Om}(\p_t\vr)^2\dx +\f{1}{2}\int_{\Om}(\Div(\vr\vu))^2\dx.
\]
Consequently,
\beq\label{pq23}
\ep \f{d}{dt}\int_{\Om}|\nabla \vr|^2 \dx +\int_{\Om}(\p_t\vr)^2\dx
\leq \int_{\Om}(\Div(\vr\vu))^2\dx.
\eeq
Notice that the classical $L^p-L^q$ estimate for parabolic equations (see \cite{AM}) indicates that
\[
\|\nabla \vr\|_{L^4(0,T;L^4(\Om))} \leq \Lambda (\|\vr \vu\|_{L^4(0,T;L^4(\Om))}+\|\vr_0\|_{H^1(\Om)}) \leq \Lambda,
\]
where (\ref{pq20}) and (\ref{pq22}) were employed. Thus, it holds that
\beq\label{pq24}
\|\Div(\vr\vu)\|_{L^2(0,T;L^2(\Om))}\leq \|\vr \Div \vu\|_{L^2(0,T;L^2(\Om))}
+\|\vu\cdot \nabla \vr\|_{L^2(0,T;L^2(\Om))}\leq \Lambda.
\eeq
In view of (\ref{pq23})-(\ref{pq24}), we finally get
\beq\label{pq25}
\|\p_t \vr\|_{L^2(0,T;L^2(\Om))} +\|\nabla \vr\|_{L^{\infty}(0,T;L^2(\Om))}\leq \Lambda.
\eeq
Writing (\ref{pq1}) as
\[
\ep \Delta \vr= \p_t \vr +\Div(\vr \vu).
\]
It follows from the classical estimates for elliptic equations and (\ref{pq24})-(\ref{pq25}) that
\beq\label{pq26}
\|\vr\|_{L^2(0,T;H^2(\Om))} \leq \Lambda.
\eeq

Following the same line as for the density, we get the estimate for the magnetic field
\beq\label{pq27}
\|\p_t b\|_{L^2(0,T;L^2(\Om))}+\|b\|_{L^2(0,T;H^2(\Om))} +\|\nabla b\|_{L^{\infty}(0,T;L^2(\Om))}\leq \Lambda.
\eeq

\noindent{\emph{ Step 6. Higher order energy estimates of the velocity.}}

By taking inner product of (\ref{pq2}) with $\p_t\vu$,
\[
\f{1}{2}\f{d}{dt} \int_{\Om} \left(  \mu|\nabla \vu|^2 +(\mu+\lambda)(\Div \vu)^2        \right)\dx
 +\int_{\Om} \vr |\p_t\vu|^2 \dx
\]
\beq\label{pq28}
=-\int_{\Om} \left( \vr \vu \cdot \nabla \vu + \ep \nabla \vr \cdot \nabla \vu +\delta\nabla (\vr+b)^{\Gamma}
+\Grad{\left(a\vr^{\gamma}+\f{1}{2}b^2\right)} \right)\cdot \p_t \vu \dx;
\eeq
whence the right-hand side can be estimated through
\[
\left|  \int_{\Om} \left( \vr \vu \cdot \nabla \vu + \ep \nabla \vr \cdot \nabla \vu +\delta\nabla (\vr+b)^{\Gamma}
+\Grad{\left(a\vr^{\gamma}+\f{1}{2}b^2\right)} \right)\cdot \p_t \vu \dx               \right|
\]
\[
\leq \nu \|\p_t \vu\| _{L^2(\Om)}^2 + \nu \|\vu\|_{H^2(\Om)}^2 + \Lambda (\nu) \left( \|\vu\|_{H^1(\Om)}^2 +\|\vr\|_{H^2(\Om)}^2 \right) \|\vu\|_{H^1(\Om)}^2
\]
\[
+ \Lambda(\nu) \left(\|\nabla \vr\| _{L^2(\Om)}^2+\|\nabla b\| _{L^2(\Om)}^2\right).
\]

To proceed, we rewrite (\ref{pq2}) as the Lam\'{e} system as follows.
\[
-\mu \Delta \vu -(\mu+\lambda)\nabla \Div \vu =
\]
\[
- \left(\p_t(\vr\vu)+\Div(\vr\vu \otimes \vu)+\Grad{\left(a\vr^{\gamma}+\f{1}{2}b^2\right)}
+\ep \nabla \vr \cdot \nabla \vu +\delta\nabla (\vr+b)^{\Gamma} \right).
\]
It then follows from the classical estimate for Lam\'{e} system (see \cite{NS}) that
\[
\|\vu\|_{H^2(\Om)}^2 \leq \Lambda \|\p_t \vu\| _{L^2(\Om)}^2
+ \Lambda   \left( \|\vu\|_{H^1(\Om)}^2 +\|\vr\|_{H^2(\Om)}^2 \right) \|\vu\|_{H^1(\Om)}^2
\]
\beq\label{pq29}
+ \Lambda \left(\|\nabla \vr\| _{L^2(\Om)}^2+\|\nabla b\| _{L^2(\Om)}^2\right).
\eeq
Finally, by inserting (\ref{pq29}) into (\ref{pq28}) and choosing $\nu$ suitably small, we deduce via a Gronwall's type argument that
\beq\label{pq30}
\|\p_t \vu\|_{L^2(0,T;L^2(\Om))} +\|\vu\|_{L^{\infty}(0,T;H^1(\Om))} + \|\vu\|_{L^{2}(0,T;H^2(\Om))} \leq \Lambda,
\eeq
where we have used (\ref{pq20}) and (\ref{pq25})-(\ref{pq27}).

\section{Passing to the limit $\ep\rightarrow 0$}\label{sen3}

\subsection{Uniform-in-$\ep$ estimates} \label{sen31}
We start our construction of weak solutions to  (\ref{p6})-(\ref{p10}) from now on. To this end, similar to \cite{FNP,NP,VWY}, we consider approximate initial data $(\vr_0^{\delta},\vu_0^{\delta},b_0^{\delta})$ satisfying
\beq\label{ps1}
0<\delta \leq \vr_0^{\delta}, b_0^{\delta} \leq \delta^{-\f{1}{2\Gamma} } ,\,\,
(\vr_0^{\delta},b_0^{\delta})\in C^3(\overline{\Om}),
\eeq
\beq\label{ps2}
\nabla \vr_0^{\delta} \cdot \vc{n}|_{\p \Om}=\nabla b_0^{\delta} \cdot \vc{n}|_{\p \Om}=0,
\eeq
\beq\label{ps3}
\vr_0^{\delta} \rightarrow \vr_0 \,\,\text{ strongly in } L^{\gamma}(\Om),
\eeq
\beq\label{ps4}
b_0^{\delta} \rightarrow b_0 \,\,\text{ strongly in } L^{2}(\Om),
\eeq
\beq\label{ps5}
\sqrt{\vr_0^{\delta}} \vu_0^{\delta} \rightarrow \sqrt{\vr_0} \vu_0 \,\,\text{ strongly in } L^{2}(\Om),\,\,\vu_0^{\delta} \in C_c^3(\Om),
\eeq
\beq\label{ps6}
C_{\star}\vr_0^{\delta} \leq b_0^{\delta} \leq C^{\star} \vr_0^{\delta},
\eeq
for $(\vr_0,\vu_0,b_0)$ subject to (\ref{p12})-(\ref{p13}) and positive constants $C_{\star},  C^{\star}$ defined in Proposition \ref{prop}. It follows from Proposition \ref{prop} that there exists a unique strong solution $(\vr_{\ep},\vu_{\ep},b_{\ep})$ to the first level approximate problem (\ref{pq1})-(\ref{pq5}) with approximate initial data $(\vr_0^{\delta},\vu_0^{\delta},b_0^{\delta})$. Moreover, as a direct consequence of the energy inequality (\ref{pq11}), we obtain the following uniform-in-$\ep$ estimates provided that $\Gamma$ is suitably large, say $\Gamma > \max\{4,\gamma\}$.
\beq\label{ps7}
\sup_{t\in(0,T)} \|\vr_{\ep}(t)\|_{L^{\gamma}(\Om)}^{\gamma} \leq C,\,\,\,
\vr_{\ep} >0 \,\text{ in } \Om\times (0,T),
\eeq
\beq\label{ps8}
\sup_{t\in(0,T)} \|b_{\ep}(t)\|_{L^{2}(\Om)}^{2} \leq C,\,\,\,
b_{\ep} >0 \,\text{ in } \Om\times (0,T),
\eeq
\beq\label{ps9}
 \delta \sup_{t\in(0,T)} \|(\vr_{\ep}(t),b_{\ep}(t))\|_{L^{\Gamma}(\Om)}^{\Gamma} \leq C,
\eeq
\beq\label{ps10}
\sup_{t\in(0,T)} \|\sqrt{\vr_{\ep}}  \vu_{\ep}(t)\|_{L^{2}(\Om)}^{2} \leq C,
\eeq
\beq\label{ps11}
\int_0^T \|\vu_{\ep}(t)\|^2_{ H_0^1(\Om) } \dt \leq C,
\eeq
\beq\label{ps12}
\ep \int_0^T \| (\nabla \vr_{\ep},\nabla b_{\ep} ) \|_{L^{2}(\Om)}^{2} \dt\leq C.
\eeq
Throughout this section, we denote by $C$ various positive constants independent of $\ep$. Notice that we conclude from (\ref{pq9}) that
\beq\label{ps13}
C_{\star}\vr_{\ep}\leq b_{\ep}\leq C^{\star}\vr_{\ep} \, \,\text{   in }\Om\times (0,T).
\eeq
Combining (\ref{ps13}) with (\ref{ps7})-(\ref{ps8}),
\beq\label{ps14}
\sup_{t\in(0,T)} \|(\vr_{\ep}(t),b_{\ep}(t))\|_{L^{ \max\{2,\gamma  \} }(\Om)}^{\max\{2,\gamma  \}} \leq C.
\eeq

In order to pass to the limit $\ep \rightarrow 0$ in the approximate equation (\ref{pq2}), we need higher integrability of $\{\vr_{\ep}\}_{\ep>0}$ and $\{b_{\ep}\}_{\ep>0}$ uniformly in $\ep$. This is accomplished by employing Bogovskii operator (see \cite{BG}). Since the proof follows the same line as \cite{FNP,VWY}, we omit the details here. Specifically, we have
\beq\label{ps15}
\int_0^T\int_{\Om} \left( \vr_{\ep}^{\gamma+1}+b_{\ep}^3 + \delta \vr_{\ep}^{\Gamma+1} +\delta b_{\ep}^{\Gamma+1}    \right)\dxdt \leq C.
\eeq

\subsection{Strong convergence of $\{\vr_{\ep}\}_{\ep>0}$ and $\{b_{\ep}\}_{\ep>0}$} \label{sen32}
With the uniform-in-$\ep$ estimates established in Section \ref{sen31}, we are now in a position to pass to the limit $\ep\rightarrow 0$ for the first level approximate problem (\ref{pq1})-(\ref{pq5}) with approximate initial data $(\vr_0^{\delta},\vu_0^{\delta},b_0^{\delta})$. Clearly, it follows from (\ref{ps7})-(\ref{ps15}) that, up to a suitable subsequence, there exists a weak limit $(\vr,\vu,b)$ such that as $\ep \rightarrow 0$:
\beq\label{ps16}
(\vr_{\ep},b_{\ep})\rightarrow (\vr,b) \text{ in } C([0,T]; L^{\Gamma}_{w}(\Om)), \text{ weakly in } L^{\Gamma+1}(\Om \times(0,T)),
\eeq
\beq\label{ps17}
(\ep \nabla \vr_{\ep},\ep \nabla b_{\ep})\rightarrow \vc{0} \,\,\text{ strongly in } L^2(0,T;L^2(\Om)),
\eeq
\beq\label{ps18}
\vu_{\ep}\rightarrow \vu \,\,\text{ weakly in }L^2(0,T;H_0^1(\Om)),
\eeq
\beq\label{ps19}
\vr_{\ep}\vu_{\ep} \rightarrow \vr \vu \,\,\text{ in } C([0,T]; L^{\f{2\gamma}{\gamma+1}  }_{w}(\Om)),
\eeq
\beq\label{ps20}
b_{\ep}\vu_{\ep}\rightarrow b \vu \,\,\text{ in }\mathcal{D}'(\Om\times (0,T)),
\eeq
\beq\label{ps21}
\vr_{\ep}\vu_{\ep}\otimes \vu_{\ep}\rightarrow \vr \vu\otimes \vu \,\,\text{ in }\mathcal{D}'(\Om\times (0,T)),
\eeq
\beq\label{ps22}
\ep  \nabla \vr_{\ep} \cdot  \nabla \vu_{\ep}\rightarrow \vc{0} \,\,\text{ strongly in } L^1(\Om \times (0,T)),
\eeq
\beq\label{ps23}
\left(a\vr_{\ep}^{\gamma}+\f{1}{2}b_{\ep}^2\right)
+\delta(\vr_{\ep}+b_{\ep})^{\Gamma} \rightarrow
\overline{ \left(a\vr^{\gamma}+\f{1}{2}b^2\right)
+\delta(\vr+b)^{\Gamma} }  \,\,\text{ weakly in } L^{\f{\Gamma+1}{\Gamma}}(\Om \times(0,T)),
\eeq
where we used the convention that $\overline{f}$ signifies the weak limit in $L^1({\Om\times (0,T)})$ of the sequence $\{f_{\ep}\}_{\ep >0}$ as $\ep \rightarrow 0$. Moreover, it follows immediately from (\ref{ps13}) and (\ref{ps16}) that
\beq\label{ps24}
\vr\geq 0, \,\,\,C_{\star}\vr\leq b\leq C^{\star}\vr\, \,\, \text{ a.e. in }\Om\times (0,T).
\eeq

As a consequence, using (\ref{ps16})-(\ref{ps23}), we pass to the limit $\ep\rightarrow 0$ in (\ref{pq1})-(\ref{pq3}) to deduce that the weak limit $(\vr,\vu,b)$ obeys
\beq\label{ps25}
\p_t\vr + \Div(\vr\vu)=0,
\eeq
\beq\label{ps26}
\p_t(\vr\vu)+\Div(\vr\vu \otimes \vu)+\Grad{\overline{ \left(a\vr^{\gamma}+\f{1}{2}b^2\right)
+\delta(\vr+b)^{\Gamma} }}=\mu\Delta\vu + (\mu +\lambda)\Grad{\Div \vu},
\eeq
\beq\label{ps27}
\p_t b+ \Div(b\vu)=0,
\eeq
in $\mathcal{D}'(\Om\times (0,T))$. In addition, the equations (\ref{ps25})-(\ref{ps27}) are supplemented with the initial and boundary conditions
\beq\label{ps28}
(\vr,\vr \vu,b)(\cdot,0)=(\vr_0^{\delta},\vr_0^{\delta}\vu_0^{\delta},b_0^{\delta}),
\eeq
\beq\label{ps29}
\vu|_{\p\Om}=\vc{0},
\eeq
in accordance with (\ref{ps16}) and (\ref{ps19}).

Obviously, it remains to show
\beq\label{ps30}
\overline{ \left(a\vr^{\gamma}+\f{1}{2}b^2\right)
+\delta(\vr+b)^{\Gamma} } = \left(a\vr^{\gamma}+\f{1}{2}b^2\right)
+\delta(\vr+b)^{\Gamma}\,\,\, \text{ a.e. in }\Om\times (0,T).
\eeq
Compared with the compressible Navier-Stokes equations of isentropic fluids, where the pressure depends on the single density, the present situation is much more complicated. When the pressure depends on the single density monotonically, pointwise convergence of approximate densities can be justified based on the nice properties enjoyed by the effective viscous flux \cite{MM,SD1} and a modified Minty's trick \cite{FNP}. We refer to \cite{FE3} for the treatment of pressure depending on the single density non-monotonically. However, in the present setting the pressure term depends on two variables with each of them satisfying the continuity equation, which makes the technique in \cite{FNP} unavailable directly. To overcome this, we appeal to the new idea of variable reduction developed by Vasseur et al. \cite{VWY} in the context of compressible two-fluid model. To this end, we first establish the following lemma indicating the nice property of the effective viscous flux.
\begin{Lemma}\label{lem31}
It holds that
\[
\lim_{\ep\rightarrow 0}\int_0^T \psi \int_{\Om} \phi \left\{\left(a\vr_{\ep}^{\gamma}+\f{1}{2}b_{\ep}^2\right)
+\delta(\vr_{\ep}+b_{\ep})^{\Gamma}-(\lambda+2\mu)\Div \vu_{\ep} \right\}
(\vr_{\ep}+b_{\ep})\dxdt
\]
\beq\label{ps31}
=\int_0^T \psi \int_{\Om} \phi \left\{ \overline{ \left(a\vr^{\gamma}+\f{1}{2}b^2\right)
+\delta(\vr+b)^{\Gamma} } -(\lambda+2\mu)\Div \vu \right\} (\vr+b)\dxdt,
\eeq
for any $\psi \in C_c^{\infty}((0,T)),\,\phi \in C_c^{\infty}(\Om)$.
\end{Lemma}

The proof of Lemma \ref{lem31} is based on the Div-Curl lemma and a suitable choice of test function in the momentum equation. For brevity, we will give the sketch of proof for Lemma \ref{lem41} in the same spirit.

With Lemma \ref{lem31} at hand, the verification of (\ref{ps30}) is carried out basically by virtue of the renormalization technique in spirit of DiPerna-Lions \cite{DL}. However, since now the pressure term depends on two variables, the classical argument for the compressible Navier-Stokes system \cite{FNP} needs to be modified. This is the main contribution of the following lemma.
\begin{Lemma}\label{lem32}
It holds that
\beq\label{ps32}
\lim_{\ep\rightarrow 0}\int_0^T\int_{\Om}(\vr_{\ep}+b_{\ep})
\left|  \f{\vr_{\ep}}{\vr_{\ep}+b_{\ep}}- \f{\vr}{\vr+b}   \right|^p \dxdt=0,
\eeq
\beq\label{ps33}
\lim_{\ep\rightarrow 0}\int_0^T\int_{\Om}(\vr_{\ep}+b_{\ep})
\left|  \f{b_{\ep}}{\vr_{\ep}+b_{\ep}}- \f{b}{\vr+b}   \right|^p \dxdt=0,
\eeq
for any $1<p<\infty$.
\end{Lemma}

The detailed proof of Lemma \ref{lem32} can be found in Theorem 2.2 of \cite{VWY}. Nevertheless, we would like to point out two key observations in the proof for the convenience of the reader. The first observation is concerned with the weak limit function $(\vr,\vu,b)$. More precisely, the limit pairs $(\vr,\vu)$ and $(b,\vu)$ satisfy the continuity equations (\ref{ps25}) and (\ref{ps27}) in the sense of distributions respectively. By defining
\[
G(\vr,b):=\f{\vr^2}{\vr+b},
\]
it follows from the generalized form of renormalization technique for transport equations (see Lemma 2.5 in \cite{VWY}) that
\[
\p_t G(\vr,b)+ \Div(G(\vr,b)\vu)+\left( \vr \p_{\vr}G(\vr,b) +b \p_{b}G(\vr,b) -G(\vr,b)  \right)\Div \vu=0
\]
in the sense of distributions. A straightforward computation gives that
\[
 \vr \p_{\vr}G(\vr,b) +b \p_{b}G(\vr,b) -G(\vr,b)=0.
\]
Hence, we find that the quantity
\beq\label{ps34}
\int_{\Om} G(\vr,b) \dx
\eeq
is independent of $t$. The second observation is related to the approximate solution $(\vr_{\ep},\vu_{\ep},b_{\ep})$. By setting
\[
F(\vr_{\ep},d_{\ep}):=\f{\vr_{\ep}^2}{d_{\ep}}=\f{\vr_{\ep}^2}{\vr_{\ep}+b_{\ep}},
\]
it follows that
\[
\p_t F(\vr_{\ep},d_{\ep})+\Div(F(\vr_{\ep},d_{\ep})\vu_{\ep})
+\left( \vr_{\ep}\p_{\vr_{\ep}} F(\vr_{\ep},d_{\ep})  +d_{\ep}\p_{d_{\ep}} F(\vr_{\ep},d_{\ep})
-  F(\vr_{\ep},d_{\ep}) \right)\Div \vu_{\ep}
\]
\[
+\ep \left( \p^2_{\vr_{\ep}^2}F(\vr_{\ep},d_{\ep}) |\nabla \vr_{\ep}|^2 +\p^2_{d_{\ep}^2}F(\vr_{\ep},d_{\ep}) |\nabla d_{\ep}|^2
+2\p^2_{\vr_{\ep}d_{\ep}   }  F(\vr_{\ep},d_{\ep}) \nabla \vr_{\ep}\cdot \nabla d_{\ep}\right)
\]
\[
-\ep \Delta  F(\vr_{\ep},d_{\ep})=0.
\]
From the definition one sees that
\[
 \vr_{\ep}\p_{\vr_{\ep}} F(\vr_{\ep},d_{\ep})  +d_{\ep}\p_{d_{\ep}} F(\vr_{\ep},d_{\ep})
-  F(\vr_{\ep},d_{\ep}) =0,
\]
which combined with the convexity of $F(\vr_{\ep},d_{\ep})$ yields that
\beq\label{ps35}
\int_{\Om} F(\vr_{\ep},d_{\ep})\dx
\eeq
is non-increasing with respect to time.

To proceed, we follow basically the argument from compressible Navier-Stokes system \cite{FNP}. By applying the renormalization technique to the continuity equations with artificial diffusion terms (\ref{pq1}) and (\ref{pq3}), it is easy to obtain the following lemma. The reader may refer to \cite{VWY} for the detailed proof.
\begin{Lemma}\label{lem33}
\[
\int_{\Om} \left( \vr_{\ep}\log \vr_{\ep} -\vr \log \vr +b_{\ep}\log b_{\ep} -b\log b \right)(\cdot,t)\dx
\]
\beq\label{ps36}
\leq \int_0^t\int_{\Om} (\vr+b)\Div \vu \dx\ds-\int_0^t\int_{\Om} (\vr_{\ep}+b_{\ep})\Div \vu_{\ep} \dx\ds,\,\, \text{ for a.e. }t\in (0,T).
\eeq
\end{Lemma}

In order to prove (\ref{ps30}), by virtue of Lemma \ref{lem33} and the convexity of function $\Phi(z)=z\log z$, it suffices to verify
\beq\label{ps37}
\int_0^t\int_{\Om} (\vr+b)\Div \vu \dx\ds\leq \lim_{\ep\rightarrow 0}\int_0^t\int_{\Om} (\vr_{\ep}+b_{\ep})\Div \vu_{\ep} \dx\ds.
\eeq
By exploiting the effective viscous flux identity formulated in Lemma \ref{lem31}, it turns out that the verification of (\ref{ps37}) reduces to (see Lemma 4.7 in \cite{VWY})
\begin{Lemma}\label{lem34}
\beq\label{ps38}
\int_0^t\psi\int_{\Om}\phi (\vr+b)\overline{ \left(a\vr^{\gamma}+\f{1}{2}b^2\right) }\dx\ds
\leq \int_0^t\psi\int_{\Om}\phi \overline{ (\vr+b)  \left(a\vr^{\gamma}+\f{1}{2}b^2\right)   }\dx\ds,
\eeq
for a.e. $t \in (0,T)$ and any nonnegative $\psi \in C_c^{\infty}((0,T)),\,\phi \in C_c^{\infty}(\Om)$.
\end{Lemma}

We remark that for the compressible isentropic Navier-Stokes equations, (\ref{ps38}) takes the form \[
\int_0^t\psi\int_{\Om}\phi \vr \overline{\vr^{\gamma}}\dx\ds
\leq \int_0^t\psi\int_{\Om}\phi \overline{\vr^{\gamma+1}}\dx\ds,
\]
which is a simple consequence of the well-known result for weak convergence and monotonicity \cite{FN1}. However, since now the pressure term depends on two variables, a careful analysis is needed to exclude possible oscillations resulting from the cross terms. The new method of variable reduction proposed in \cite{VWY} is based on the following decomposition:
\[
a \vr_{\ep}^{\gamma}+\f{1}{2}b_{\ep}^2= a A_{\ep} ^{\gamma} d_{\ep}^{\gamma}+\f{1}{2}B_{\ep}^2 d_{\ep}^2=a A^{\gamma} d_{\ep}^{\gamma}+\f{1}{2}B^2 d_{\ep}^2
+a(A_{\ep} ^{\gamma}-A^{\gamma})d_{\ep}^{\gamma} +\f{1}{2}(B_{\ep}^2 -B^2 )d_{\ep}^2,
\]
\[
\vr_{\ep}+b_{\ep}=(A_{\ep}+B_{\ep})d_{\ep} =(A+B)d_{\ep}+(A_{\ep}+B_{\ep}-A-B)d_{\ep},
\]
where $d_{\ep}:=\vr_{\ep}+b_{\ep},\,d:=\vr+b,\,A_{\ep}:=\f{\vr_{\ep}}{d_{\ep}},\,B_{\ep}:=\f{b_{\ep}}{d_{\ep}},
\,A:=\f{\vr}{d},\,B:=\f{b}{d}$. Then, by Lemma \ref{lem32} and the uniform-in-$\ep$ estimates obtained in Section \ref{sen31},
\[
\int_0^t\psi\int_{\Om}\phi \left(a A^{\gamma} d_{\ep}^{\gamma}+\f{1}{2}B^2 d_{\ep}^2\right)(A_{\ep}+B_{\ep}-A-B)d_{\ep} \dx\ds
\]
and
\[
\int_0^t\psi\int_{\Om}\phi \left( a(A_{\ep} ^{\gamma}-A^{\gamma})d_{\ep}^{\gamma} +\f{1}{2}(B_{\ep}^2 -B^2 )d_{\ep}^2 \right)(\vr_{\ep}+b_{\ep}) \dx\ds
\]
vanish as $\ep \rightarrow 0$ for any $\psi \in C_c^{\infty}((0,T)),\,\phi \in C_c^{\infty}(\Om)$. We refer to Lemma 4.5 in \cite{VWY} for the details.

\section{Passing to the limit $\delta\rightarrow 0$}\label{sen4}

\subsection{Uniform-in-$\delta$ estimates} \label{sen41}
We conclude from Section \ref{sen3} that, upon passing to the limit $\ep\rightarrow 0$, there exists a weak solution $(\vr_{\delta},\vu_{\delta},b_{\delta})$ satisfying
\beq\label{pt1}
\p_t\vr + \Div(\vr\vu)=0,
\eeq
\beq\label{pt2}
\p_t(\vr\vu)+\Div(\vr\vu \otimes \vu)+\Grad{ \left(a\vr^{\gamma}+\f{1}{2}b^2\right)
+\delta\nabla(\vr+b)^{\Gamma} }=\mu\Delta\vu + (\mu +\lambda)\Grad{\Div \vu},
\eeq
\beq\label{pt3}
\p_t b+ \Div(b\vu)=0,
\eeq
in $\mathcal{D}'(\Om\times (0,T))$, together with the initial and boundary conditions:
\beq\label{pt4}
(\vr,\vr \vu,b)(\cdot,0)=(\vr_0^{\delta},\vr_0^{\delta}\vu_0^{\delta},b_0^{\delta}),
\eeq
\beq\label{pt5}
\vu|_{\p\Om}=\vc{0}.
\eeq
Thus, our ultimate goal is to pass to the limit $\delta\rightarrow 0$ to recover a weak solution to the original problem (\ref{p6})-(\ref{p10}). To this end, we first notice that the following uniform-in-$\delta$ estimates are direct consequences of Section \ref{sen31} and the weak convergence results in Section \ref{sen32}. In this section, various positive constants are represented by the same letter $C$ independent of $\delta$.
\beq\label{pt6}
\sup_{t\in(0,T)} \|\vr_{\delta}(t)\|_{L^{\gamma}(\Om)}^{\gamma} \leq C,\,\,\,
\vr_{\delta} \geq 0 \,\text{ in } \Om\times (0,T),
\eeq
\beq\label{pt7}
\sup_{t\in(0,T)} \|b_{\delta}(t)\|_{L^{2}(\Om)}^{2} \leq C,\,\,\,
b_{\delta} \geq 0 \,\text{ in } \Om\times (0,T),
\eeq
\beq\label{pt8}
 \delta \sup_{t\in(0,T)} \|(\vr_{\delta}(t),b_{\delta}(t))\|_{L^{\Gamma}(\Om)}^{\Gamma} \leq C,
\eeq
\beq\label{pt9}
\sup_{t\in(0,T)} \|\sqrt{\vr_{\delta}}  \vu_{\delta}(t)\|_{L^{2}(\Om)}^{2} \leq C,
\eeq
\beq\label{pt10}
\int_0^T \|\vu_{\delta}(t)\|^2_{ H_0^1(\Om) } \dt \leq C,
\eeq
\beq\label{pt11}
C_{\star}\vr_{\delta}\leq b_{\delta}\leq C^{\star}\vr_{\delta}\, \, \text{   a.e. in }\Om\times (0,T),
\eeq
\beq\label{pt12}
\sup_{t\in(0,T)} \|(\vr_{\delta}(t),b_{\delta}(t))\|_{L^{ \max\{2,\gamma  \} }(\Om)}^{\max\{2,\gamma  \}} \leq C.
\eeq
Similar to (\ref{ps15}), we need higher integrability of $\{\vr_{\delta}\}_{\delta>0}$ and $\{b_{\delta}\}_{\delta>0}$ uniformly in $\delta$ so as to pass to the limit $\delta \rightarrow 0$ in (\ref{pt2}). To do this, again by virtue of Bogovskii operator and the uniform-in-$\delta$ estimates (\ref{pt6})-(\ref{pt12}), we have
\beq\label{pt13}
\int_0^T\int_{\Om}
\left(  \vr_{\delta}^{\gamma+\vt} +b_{\delta}^{2+\vt}+\delta \vr_{\delta}^{\Gamma+\vt} +\delta b_{\delta}^{\Gamma+\vt} \right)
 \dxdt \leq C,
\eeq
for some $\vt>0$ depending only on $\gamma$ (see \cite{FNP}). In light of (\ref{pt11}), it holds that
\beq\label{pt15}
\int_0^T\int_{\Om}  \left(  \vr_{\delta}^{  \max\{\gamma+\vt,2+\vt  \} }+b_{\delta}^{\max\{\gamma+\vt,2+\vt  \} }        \right)\dxdt
\leq C.
\eeq
With the above uniform-in-$\delta$ estimates at hand, we can pass to the limit $\delta\rightarrow 0$ in the approximate system (\ref{pt1})-(\ref{pt3}) to conclude, up to a suitable subsequence, that
\beq\label{pt16}
\vr_{\delta}\rightarrow \vr \text{ in }C([0,T];L_w^{ \max\{\gamma,2\}   }(\Om)),\,\,
\text{ weakly in }L^{ \max\{\gamma+\vt,2+\vt \} }(\Om \times(0,T)),
\eeq
\beq\label{pt17}
b_{\delta}\rightarrow b \text{ in }C([0,T];L_w^{ \max\{\gamma,2\}   }(\Om)),\,\,
\text{ weakly in }L^{ \max\{\gamma+\vt,2+\vt \} }(\Om \times(0,T)),
\eeq
\beq\label{pt18}
\vu_{\delta}\rightarrow \vu \text{  weakly in }L^2(0,T;H_0^1(\Om)),
\eeq
\beq\label{pt19}
\vr_{\delta}\vu_{\delta} \rightarrow \vr \vu \,\,\text{ in } C([0,T]; L^{\f{2\max\{\gamma,2\} }{\max\{\gamma,2\} +1}  }_{w}(\Om)),
\eeq
\beq\label{pt20}
b_{\delta}\vu_{\delta}\rightarrow b \vu \,\,\text{ in }\mathcal{D}'(\Om\times (0,T)),
\eeq
\beq\label{pt21}
\vr_{\delta}\vu_{\delta}\otimes \vu_{\delta}\rightarrow \vr \vu\otimes \vu \,\,\text{ in }\mathcal{D}'(\Om\times (0,T)),
\eeq
\beq\label{pt22}
\delta (\vr_{\delta}+b_{\delta})^{\Gamma}\rightarrow 0 \,\text{  strongly in }L^1(\Om\times(0,T)),
\eeq
\beq\label{pt23}
\vr_{\delta}^{\gamma}\rightarrow \overline{\vr^{\gamma}}\, \,\text{ weakly in } L^{ \f{\max\{\gamma+\vt,2+\vt  \}}  {\gamma}   }(\Om \times(0,T)),
\eeq
\beq\label{pt24}
b_{\delta}^2\rightarrow \overline{ b^2 }\,\,\text{ weakly in } L^ {    \f{\max\{\gamma+\vt,2+\vt  \}}{2}     }(\Om\times(0,T)).
\eeq

Therefore, we realize that the weak limit function $(\vr,\vu,b)$ obeys
\beq\label{pt25}
\p_t\vr + \Div(\vr\vu)=0,
\eeq
\beq\label{pt26}
\p_t(\vr\vu)+\Div(\vr\vu \otimes \vu)+\Grad{\left( a\overline{\vr^{\gamma}}+\f{1}{2}\overline{b^2}\right)}=\mu\Delta\vu + (\mu +\lambda)\Grad{\Div \vu},
\eeq
\beq\label{pt27}
\p_t b+ \Div(b\vu)=0,
\eeq
in $\mathcal{D}'(\Om\times(0,T))$ supplemented with the initial and boundary conditions:
\beq\label{pt28}
(\vr,\vr \vu,b)(\cdot,0)=(\vr_0,\vr_0\vu_0,b_0),
\eeq
\beq\label{pt29}
\vu|_{\p\Om}=\vc{0}.
\eeq
Obviously, to finish the proof of Theorem \ref{ls}, it remains to verify
\beq\label{pt30}
\overline{\vr^{\gamma}}=\vr^{\gamma} \,\,\text{ a.e. in } \Om\times(0,T),
\eeq
\beq\label{pt31}
 \overline{ b^2 }=b^2 \,\,\text{ a.e. in } \Om\times(0,T),
\eeq
which is the goal of the next subsection.

\subsection{Strong convergence of $\{\vr_{\delta}\}_{\delta>0}$ and $\{b_{\delta}\}_{\delta>0}$}\label{sen42}

The main idea to show pointwise convergence of $\{\vr_{\delta}\}_{\delta>0}$ and $\{b_{\delta}\}_{\delta>0}$ is the same as the limit passage for $\ep\rightarrow 0$. As in Section \ref{sen32}, the key issue is uniform-in-$\ep$ estimates of $\{\vr_{\ep}\}_{\ep>0}$ and $\{b_{\ep}\}_{\ep>0}$ in $L^2(0,T;L^2(\Om))$ in order to apply the generalized form of DiPerna-Lions' argument for the renormalized solutions to transport equations (see Lemma 2.5 in \cite{VWY}). It turns out that this condition is naturally satisfied due to (\ref{pt15}). Remarkably, for the compressible two-fluid model discussed in \cite{VWY}, the adiabatic exponents are chosen in such a way that the two densities are uniformly bounded in $L^2(0,T;L^2(\Om))$, while in the classical work \cite{FNP}, the boundedness of oscillations is investigated so as to verify that the limit pair $(\vr,\vu)$ obeys the continuity equation in the renormalized sense.

For the sake of simplicity, we only point out the difference between the two limit passages. To begin with, we notice that the nice property admitted by the effective viscous flux still plays a crucial role at this stage. However, due to the low integrability of $\{\vr_{\delta}^{\gamma}\}_{\delta>0}$ and $\{b_{\delta}^{2}\}_{\delta>0}$, we need to make use of the special cut-off functions introduced in \cite{FE1}, i.e.,
\[
T_k(z):=kT\left(\f{z}{k}\right),\,\,k\geq 1,
\]
where $T(\cdot)$ is a smooth concave function on $[0,\infty)$ satisfying
\begin{equation*}
T(z)=
\left\{
\begin{aligned}
z,\text{ if }z\in [0,1], \\
2,\text{ if }z\in [3,\infty). \\
\end{aligned}
\right.
\end{equation*}

With the cut-off functions $T_k(z)$ defined above, we have the following lemma indicating the nice property of the effective viscous flux. Notice that the first two terms in the momentum equation are different from the compressible two-fluid model \cite{VWY}. Hence the analysis here is more involved.
\begin{Lemma}\label{lem41}
It holds that
\[
\lim_{\delta\rightarrow 0}\int_0^T \psi \int_{\Om} \phi \left\{\left(a\vr_{\delta}^{\gamma}+\f{1}{2}b_{\delta}^2\right)
-(\lambda+2\mu)\Div \vu_{\delta} \right\}
\left(T_k(\vr_{\delta})+T_k(b_{\delta})\right)\dxdt
\]
\beq\label{pt32}
=\int_0^T \psi \int_{\Om} \phi \left\{ \left( a\overline{\vr^{\gamma}}+\f{1}{2}\overline{b^2}\right) -(\lambda+2\mu)\Div \vu \right\} \left( \overline{T_k(\vr)  }+\overline{ T_k(b) } \right)\dxdt,
\eeq
for any $\psi \in C_c^{\infty}((0,T)),\,\phi \in C_c^{\infty}(\Om)$.
\end{Lemma}
{\bf Proof.} Obviously, it suffices to show
\[
\lim_{\delta\rightarrow 0}\int_0^T \psi \int_{\Om} \phi \left\{\left(a\vr_{\delta}^{\gamma}+\f{1}{2}b_{\delta}^2\right)
-(\lambda+2\mu)\Div \vu_{\delta} \right\}
T_k(\vr_{\delta})\dxdt
\]
\beq\label{pt34}
=\int_0^T \psi \int_{\Om} \phi \left\{ \left( a\overline{\vr^{\gamma}}+\f{1}{2}\overline{b^2}\right) -(\lambda+2\mu)\Div \vu \right\}  \overline{T_k(\vr)  } \dxdt,
\eeq
and
\[
\lim_{\delta\rightarrow 0}\int_0^T \psi \int_{\Om} \phi \left\{\left(a\vr_{\delta}^{\gamma}+\f{1}{2}b_{\delta}^2\right)
-(\lambda+2\mu)\Div \vu_{\delta} \right\}
T_k(b_{\delta})\dxdt
\]
\beq\label{pt35}
=\int_0^T \psi \int_{\Om} \phi \left\{ \left( a\overline{\vr^{\gamma}}+\f{1}{2}\overline{b^2}\right) -(\lambda+2\mu)\Div \vu \right\}\overline{ T_k(b) } \\dxdt.
\eeq
We focus attention on the verification of (\ref{pt35}). It follows from (\ref{pt7}) and (\ref{pt10}) that $(b_{\delta},\vu_{\delta})$ solves (\ref{pt3}) in the renormalized sense. In particular,
\beq\label{pt36}
\p_t (T_k(b_{\delta}))+\Div (T_k(b_{\delta}) \vu_{\delta})+
(T'_k(b_{\delta})b_{\delta}-T_k(b_{\delta}))\Div \vu_{\delta}=0
\eeq
in $\mathcal{D}'(\Om\times (0,T))$. Upon passing to the weak limit,
\beq\label{pt37}
\p_t (\overline{T_k(b)})+\Div (\overline{T_k(b)}\vu)+
\overline{(T'_k(b)b-T_k(b))\Div \vu}=0
\eeq
in $\mathcal{D}'(\Om\times (0,T))$. To proceed, we introduce several operators as follows. Let $\mathcal{A}_j$ be the operator $\Delta^{-1}\p_{x_j}$ corresponding to its Fourier symbol $\f{-i\xi_j}{|\xi|^2}$; $\mathcal{R}_{i,j}$ signifies the operator $\p_{x_j}\mathcal{A}_i$ corresponding to its Fourier symbol $\f{\xi_i\xi_j}{|\xi|^2}$. For simplicity, we set $\mathcal{A}:=(\mathcal{A}_1,\mathcal{A}_2)$ and $\mathcal{R}:=\{ \mathcal{R}_{i,j} \}_{i,j=1}^2$. Testing the momentum equation (\ref{pt2}) by
\[
\psi\phi\mathcal{A}[1_{\Om}T_k(b_{\delta})],
\]
we obtain after a straightforward manipulation that
\[
\int_0^T\int_{\Om} \psi\phi\left\{ \left(a\vr_{\delta}^{\gamma}+\f{1}{2}b_{\delta}^2\right)
+\delta (\vr_{\delta}+b_{\delta})^{\Gamma}-(\lambda+2\mu)\Div \vu_{\delta}\right\}T_k(b_{\delta})\dxdt
\]
\[
=\int_0^T\int_{\Om} \psi \left\{ (\lambda+\mu)\Div \vu_{\delta} - \left(a\vr_{\delta}^{\gamma}+\f{1}{2}b_{\delta}^2\right)- \delta (\vr_{\delta}+b_{\delta})^{\Gamma} \right\} \Grad \phi \cdot\mathcal{A}[1_{\Om}T_k(b_{\delta})]\dxdt
\]
\[
+\mu \int_0^T\int_{\Om} \psi \left\{  \Grad \vu_{\delta}\cdot \Grad \phi\cdot \mathcal{A}[1_{\Om}T_k(b_{\delta})]-\sum_{i,j=1}^2\left( u^i_{\delta}\p_{x_j}\phi \mathcal{R}_{i,j}[1_{\Om}T_k(b_{\delta})]\right) +\vu_{\delta}\cdot \Grad \phi T_k(b_{\delta})  \right\}\dxdt
\]
\[
-\int_0^T\int_{\Om} \phi \vr_{\delta}\vu_{\delta} \cdot\left\{ \p_t \psi \mathcal{A}[1_{\Om}T_k(b_{\delta})]+\psi \mathcal{A}[1_{\Om}(T_k(b_{\delta})-T'_k(b_{\delta})b_{\delta}) \Div \vu_{\delta}]  \right\}\dxdt
\]
\[
-\int_0^T\int_{\Om} \psi \sum_{i,j=1}^2 \vr_{\delta}u^i_{\delta} u^j_{\delta}\p_{x_j}\phi \mathcal{A}_i[1_{\Om}T_k(b_{\delta})] \dxdt
\]
\beq\label{pt38}
+\int_0^T\int_{\Om} \psi\phi\left\{ \vr_{\delta}\vu_{\delta}\cdot\mathcal{A}[1_{\Om}\Div (T_k(b_{\delta})\vu_{\delta})] -\vr_{\delta}(\vu_{\delta}\otimes \vu_{\delta}):\mathcal{R}[1_{\Om}T_k(b_{\delta})]    \right\}\dxdt.
\eeq
In a similar manner, we use
\[
\psi\phi\mathcal{A}[1_{\Om}\overline{T_k(b)}]
\]
as a test function in the limit momentum equation (\ref{pt26}) to arrive at
\[
\int_0^T\int_{\Om}\psi\phi \left\{   \left( a\overline{\vr^{\gamma}}+\f{1}{2}\overline{b^2}\right)-(\lambda+2\mu)\Div \vu \right\} \overline{T_k(b)} \dxdt
\]
\[
=\int_0^T\int_{\Om} \psi \left\{  (\lambda+\mu)\Div \vu - \left( a\overline{\vr^{\gamma}}+\f{1}{2}\overline{b^2}\right)   \right\}\Grad \phi \cdot \mathcal{A}[1_{\Om}\overline{T_k(b)}]\dxdt
\]
\[
+\mu \int_0^T\int_{\Om} \psi \left\{  \Grad \vu\cdot \Grad \phi\cdot \mathcal{A}[1_{\Om}\overline{T_k(b)} ]-\sum_{i,j=1}^2\left( u^i\p_{x_j}\phi \mathcal{R}_{i,j}[1_{\Om}\overline{T_k(b)} ]\right) +\vu\cdot \Grad \phi \overline{T_k(b)}   \right\}\dxdt
\]
\[
-\int_0^T\int_{\Om} \phi \vr\vu \cdot\left\{ \p_t \psi \mathcal{A}[1_{\Om}\overline{T_k(b)}]+\psi \mathcal{A}[1_{\Om}\overline{(T_k(b)-T'_k(b)b) \Div \vu}]  \right\}\dxdt
\]
\[
-\int_0^T\int_{\Om} \psi \sum_{i,j=1}^2 \vr u^i u^j\p_{x_j}\phi \mathcal{A}_i[1_{\Om}\overline{T_k(b)}] \dxdt
\]
\beq\label{pt39}
+\int_0^T\int_{\Om} \psi\phi\left\{ \vr\vu\cdot\mathcal{A}[1_{\Om}\Div (\overline{T_k(b)}\vu)] -\vr(\vu\otimes \vu):\mathcal{R}[1_{\Om}\overline{T_k(b)}]    \right\}\dxdt.
\eeq
Notice that the two integrals in (\ref{pt38}) containing $\delta (\vr_{\delta}+b_{\delta})^{\Gamma}$ vanish as $\delta\rightarrow 0$ due to (\ref{pt22}). Consequently, it suffices to prove that each term on the right-hand side of (\ref{pt38}) converges to its counterpart in (\ref{pt39}). It is easy to see that as $\delta\rightarrow 0$
\beq\label{pt40}
T_k(b_{\delta})\rightarrow \overline{T_k(b)} \text{ in } C([0,T];L^p_w(\Om)),
\eeq
for any $1\leq p<\infty$. By H\"{o}rmander-Mikhlin theorem and Sobolev's embedding theorem, it holds that as $\delta\rightarrow 0$
\beq\label{pt41}
\mathcal{A}[1_{\Om}T_k(b_{\delta})]\rightarrow \mathcal{A}[1_{\Om}\overline{T_k(b)}] \text{ in }C(\overline{\Om}\times [0,T]).
\eeq
Therefore, the four terms on the right-hand side of (\ref{pt38}) converge to their counterparts in (\ref{pt39}) by virtue of (\ref{pt16})-(\ref{pt24}) and (\ref{pt41}). The last term is much more involved. To this end, integration by parts yields
\[
\int_0^T\int_{\Om} \psi\phi\left\{ \vr_{\delta}\vu_{\delta}\cdot\mathcal{A}[1_{\Om}\Div (T_k(b_{\delta})\vu_{\delta})] -\vr_{\delta}(\vu_{\delta}\otimes \vu_{\delta}):\mathcal{R}[1_{\Om}T_k(b_{\delta})]    \right\}\dxdt
\]
\[
=\int_0^T\int_{\Om} \psi\phi \vu_{\delta}\cdot \left\{ T_k(b_{\delta})  \mathcal{A}[1_{\Om}\Div (\vr_{\delta}\vu_{\delta})] -\vr_{\delta}\vu_{\delta}\cdot \mathcal{R}[1_{\Om}T_k(b_{\delta})  \right\}   \dxdt
\]
\beq\label{pt42}
+\int_0^T\psi\int_{\Om}\sum_{i,j=1}^2 T_k(b_{\delta})u^j_{\delta}
\Delta^{-1}\left\{ \p_{x_j}(\vr_{\delta}u^i_{\delta})\p_{x_i}\phi
+(\p_{x_i x_j}\phi) \vr_{\delta}u^i_{\delta} + \p_{x_j}\phi \p_{x_i}(\vr_{\delta}u^i_{\delta})   \right\}  \dxdt.
\eeq
As above, the second term on the right-hand side of (\ref{pt42}) converges to its counterpart. Remarkably, the convergence of the first term on the right-hand side of (\ref{pt42}) is verified by means of the Div-Curl lemma; see \cite{FNP} for the details. Now (\ref{pt35}) follows from (\ref{pt38}) and (\ref{pt39}) by vanishing $\delta$. Finally, we point out that the proof of (\ref{pt34}) is carried out exactly as (\ref{pt35}). The details are omitted.  $\Box$

Based on Lemma \ref{lem41} and Lemma \ref{lem32}, it follows that
\begin{Lemma}\label{lem42}
\[
\int_0^t\psi\int_{\Om}\phi \left( \overline{T_k(\vr)  }+\overline{ T_k(b) } \right)\overline{ \left(a\vr^{\gamma}+\f{1}{2}b^2\right) }\dx\ds
\]
\beq\label{pt33}
\leq \int_0^t\psi\int_{\Om}\phi\, \overline{ (T_k(\vr)+T_k(b))  \left(a\vr^{\gamma}+\f{1}{2}b^2\right)   }\dx\ds,
\eeq
for a.e. $t \in (0,T)$ and any nonnegative $\psi \in C_c^{\infty}((0,T)),\,\phi \in C_c^{\infty}(\Om)$.
\end{Lemma}

Other considerations remain basically unchanged compared with Section \ref{sen32}. The reader may refer to \cite{VWY} for the details. As a consequence, we conclude that the weak limit function $(\vr,\vu,b)$ is a global weak solution to the initial-boundary value problem (\ref{p6})-(\ref{p10}), thus finishing the proof of Theorem \ref{ls}.

\centerline{\bf Acknowledgement}
The research of Y. Sun is supported by NSF of China under grant number 11571167, 11771395, 11771206  and PAPD of Jiangsu Higher Education Institutions. The research of Y. Li is partially supported by  Postgraduate Research and Practice Innovation Program of Jiangsu Province under grant number KYCX 18-0028 and China Scholarship Council under grant number 201806190103. The authors would like to thank E. Feireisl and A. Novotn\'{y} for valuable discussions. Finally, the authors are indebted to the anonymous reviewers for careful reading of the manuscript and providing very helpful comments.


\end{document}